\documentclass[10pt]{amsart}

\usepackage{amsthm}
\newtheorem*{15thm}{15-Theorem}
\newtheorem*{15thmHerm}{15-Theorem for Hermitian Lattices}
\newtheorem*{290thm}{290-Theorem}
\newtheorem*{classThm}{Main Theorem}

\newcommand\inprod[2]{\left<#1,#2\right>}
\newcommand\Q{\mathbb{Q}}
\newcommand\Qs[1]{\Q(\sqrt{#1})}

\usepackage{amssymb}
\usepackage{url}
\usepackage{amsmath}
\usepackage{amsfonts}
\usepackage{verbatim}

\title{On Universal Binary Hermitian Forms}
\author[Scott Duke Kominers]{Scott Duke Kominers}
\address{Department of Mathematics, Harvard University\newline \indent c/o 8520 Burning Tree Road, Bethesda, MD 20817}\email{kominers@fas.harvard.edu}
\subjclass[2000]{11E39 (11E20, 11E41, 11E25)}
\keywords{universal form, universality criterion, Hermitian form}
\date{\today}

\begin{document}
\maketitle

\section{Introduction}
The question of representing integers by quadratic forms dates back to the time of Fermat, whose \emph{Two Squares Theorem} solved the question of which primes could be represented by the form $x^2+y^2$ (see \cite[p. 219]{HW}).  This theorem was later generalized by Lagrange, who showed in his \textit{Four Squares Theorem} \cite{Lagrange:four} that every positive-integer can be written as a sum of four squares of integers.

Lagrange's theorem has led to the modern study of \textit{universal forms}, those forms which represent all positive integers.  In the first half of the twentieth century, Ramanujan \cite{Ra} identified the universal positive-definite classically-integral quaternary diagonal quadratic forms, up to equivalence.   Maass \cite{Ma} and Chan, Kim, and Raghavan \cite{CKR} gave analogous classification results leading to the full classification of the positive-definite classically integral ternary quadratic forms which are universal over real quadratic fields.

Motivated by the work on universal quadratic forms over real fields, Earnest and Khosravani \cite{EK} sought an analogous classification of universal binary Hermitian forms over imaginary quadratic fields.  Recently, Iwabuchi \cite{Iw} and Kim and Park~\cite{KP} finished Earnest and Khosravani's program, completing the list of universal binary Hermitian forms.

A different direction of recent research has focused on the search for \emph{universality criteria}, simple tests which characterize the universality of positive-definite quadratic forms.  The earliest-discovered result in this vein is Conway and Schneeberger's surprising \emph{15-Theorem} (see \cite{Conway:universality,Bhargava:Fif}): \begin{15thm}A positive-definite classically integral quadratic form is universal if and only if it represents the nine ``critical numbers'' $$\{1,2,3,5,6,7,10,14,15\}.$$\end{15thm}  More recently, Bhargava and Hanke (see \cite{290}) showed an analogous criterion for the universality of positive-definite nonclassically integral quadratic forms: \begin{290thm}A positive-definite nonclassically integral quadratic form is universal if and only if it represents the numbers \begin{align*}S_{290}:=\{&1, 2, 3, 5, 6, 7, 10, 13, 14, 15, 17, 19, 21, 22, 23, 26, \\&29,30, 31, 
34, 35, 37, 42, 58, 93, 110, 145, 203, 290\}.\end{align*}\end{290thm}

While the criterion theorems reduce testing a form's universality to a simple compuation, they have rarely been used in practice.   The reason for this somewhat curious fact is that the proofs of both the 15- and 290-Theorems rely on independent identification of many universal forms of low rank, called the \emph{universal escalators}.  

The results on Hermitian forms, however, give us a chance to greatly simplify prior work through an application of the 290-Theorem.  Specifically, we apply the 290-Theorem to reduce the most difficult universality verifications in Earnest and Khosravani's \cite{EK} and Iwabuchi's \cite{Iw} papers to simple, finite computations.

\section{Preliminaries}
We let $E$ be an imaginary quadratic field over $\Q$ and let $m>0$ be a squarefree integer for which $E=\Qs{-m}$.  We denote the $\Q$-involution of $E$ by $\overline{\phantom{a}}$ and the ring of integers of $E$ by $\mathcal{O}_E$.

We let $V/E$ be an $n$-dimensional Hermitian space over $E$ with nondegenerate Hermitian form $H$.  As shown by Jacobson \cite{Ja}, we may consider $(V,H)$ as a $2n$-dimensional quadratic space $(\widetilde{V},B)$ with the bilinear form $B$  defined by the trace map $$B(x,y)=\frac{1}{2}\mathrm{Tr}_{E/\Q}(H(x,y)).$$

An  $\mathcal{O}_E$-lattice $L$ is a finitely generated $\mathcal{O}_E$-module on the Hermitian space $(V,H)$.  We consider only positive-definite integral $\mathcal{O}_E$-lattices $L$, that is, those for which $H(x,y)\in\mathcal{O}_E$ for all $x,y\in L$ and $H(x,x)>0$ for all $L\ni x \neq 0$.  If an $\mathcal{O}_E$-lattice $L$ is of the form $L= L_1\oplus L_2$ for sublattices $L_1$, $L_2$ of $L$ and $\inprod{L_1}{L_2}=0$ then we write $L\cong L_1\bot L_2$.

When $E$ has class number $1$, the ring $\mathcal{O}_E$ is a principal ideal domain whereby every $\mathcal{O}_E$-lattice $L$ is free.  Hence, in this case we may think of the Hermitian form $H$ acting on $L$ as a  function $f:\mathcal{O}_E^n\to \mathbb{Z}$ defined by $$f(x_1,\ldots, x_n)=H\left(\sum_{i=1}^n x_iv_i\right)=\sum_{i=1}^n H(v_i,v_j)x_i\bar{x}_j$$ for some suitable basis $\{v_i\}_{i=1}^n$ of $L$.  If the basis $\{v_i\}_{i=1}^n$ is orthogonal, we write $L\cong \left<H(v_1),\ldots,H(v_n)\right>$.  (For example, the form $x\bar{x}+2y\bar{y}$ is associated to the lattice $\left<1,2\right>$.)

Similarly, we may associate a quadratic lattice $\widetilde{L}$ with every Hermitian $\mathcal{O}_E$-lattice $L$.  The ring $\mathcal{O}_E$ has a basis $\{1,\omega_m\}$ as a $\mathbb{Z}$-module, where $$\omega_m=\left\{\begin{array}{cc}\frac{1+\sqrt{-m}}{2} & m\equiv 3\bmod 4\\\sqrt{-m}& \mathrm{otherwise}\end{array}\right..$$  Then, $\tilde{f}(x_1,x_2,\ldots, x_n,y_n)=f(x_1+\omega_m y_1,\ldots, x_n +\omega_m y_n)$ is a quadratic form in $2n$ variables corresponding to the lattice $\widetilde{L}$.  From this construction, it is clear that the Hermitian form $f$ is universal if and only if the quadratic form $\tilde{f}$ is.  We write $\sim$ to denote the correspondence between a Hermitian lattice and its associated quadratic form.  

\section{Classification of Universal Hermitian Forms}

Earnest and Khosravani \cite{EK}, Iwabuchi \cite{Iw}, and Kim and Park \cite{KP} identified all potentially universal Hermitian forms over imaginary quadratic fields $E$. This ``screening process'' is the more straightforward part of the classification, relying on a uniform computational method (see \cite{EK}). 

The universality of the candidates identified was then shown by a variety of methods.  Indeed, a total of eight different approaches were used.  Six of these methods were ``\emph{ad hoc}'' arguments, each an intricate method developed to prove the universality of an individual Hermitian form.

We give a unified proof of the universalities of the forms in the classification, relying primarily on Ramanujan's list of universal forms~\cite{Ra} and the 290-Theorem~\cite{290}.  The universalities of twenty-four of the twenty-five universal binary Hermitian forms follow directly from our methods, whence we require only one \emph{ad hoc} argument in our proof of the classification.
\begin{classThm}
Up to equivalence, the integral positive-definite universal binary Hermitian lattices in imaginary quadratic fields are exactly the lattices in (\ref{lat}):
\begin{equation}\label{lat}
\begin{tabular}{c|l}
$\Qs{-m}$&universal binary lattices\\
\hline
$\Qs{-1}$&$\left<1,1\right>,\left<1,2\right>,\left<1,3\right>$,\\
$\Qs{-2}$&$\left<1,1\right>,\left<1,2\right>,\left<1,3\right>,\left<1,4\right>,\left<1,5\right>$,\\
$\Qs{-3}$&$\left<1,1\right>,\left<1,2\right>$,\\
$\Qs{-5}$&$\left<1,2\right>,\left<1\right>\bot\left(\begin{array}{cc}2&-1+\omega_5\\-1+\bar{\omega}_5&3\end{array}\right)$,\\
$\Qs{-6}$&$\left<1\right>\bot\left(\begin{array}{cc}2&\omega_6\\\bar{\omega}_6&3\end{array}\right)$,\\
$\Qs{-7}$&$\left<1,1\right>,\left<1,2\right>,\left<1,3\right>$,\\
$\Qs{-10}$&$\left<1\right>\bot\left(\begin{array}{cc}2&\omega_{10}\\\bar{\omega}_{10}&5\end{array}\right)$,\\
$\Qs{-11}$&$\left<1,1\right>,\left<1,2\right>$,\\ 
$\Qs{-15}$&$\left<1\right>\bot\left(\begin{array}{cc}2&\omega_{15}\\\bar{\omega}_{15}&2\end{array}\right)$,\\
$\Qs{-19}$&$\left<1,2\right>$,\\
$\Qs{-23}$&$\left<1\right>\bot\left(\begin{array}{cc}2&\omega_{23}\\\bar{\omega}_{23}&3\end{array}\right),\left<1\right>\bot\left(\begin{array}{cc}2&-1+\omega_{23}\\-1+\bar{\omega}_{23}&3\end{array}\right)$,\\
$\Qs{-31}$&$\left<1\right>\bot\left(\begin{array}{cc}2&\omega_{31}\\\bar{\omega}_{31}&4\end{array}\right),\left<1\right>\bot\left(\begin{array}{cc}2&-1+\omega_{31}\\-1+\bar{\omega}_{31}&4\end{array}\right)$.
\end{tabular}
\end{equation}
\end{classThm}
\begin{proof}[Proof of the Main Theorem]
Earnest and Khosravani \cite{EK}, Iwabuchi \cite{Iw}, and Kim and Park \cite{KP} showed that no binary Hermitian forms not in the list (\ref{lat}) can be universal over an imaginary quadratic field $E$.  Therefore, we must only show the universality of each of these candidate forms.

First, we identify the lattices in the list (\ref{lat}) which correspond to diagonal quaternary quadratic forms:
\begin{equation}\label{diags}
\begin{tabular}{rcl}
$\left<1,1\right>$ in $\Qs{-1}$ &$\sim$& $w^2+x^2+y^2+z^2$,\\
$\left<1,1\right>$ in $\Qs{-2}$ &$\sim$& $w^2+x^2+2y^2+2z^2$,\\
$\left<1,2\right>$ in $\Qs{-1}$ &$\sim$& $w^2+x^2+2y^2+2z^2$,\\
$\left<1,2\right>$ in $\Qs{-2}$ &$\sim$& $w^2+2x^2+2y^2+4z^2$,\\
$\left<1,2\right>$ in $\Qs{-5}$ &$\sim$& $w^2+2x^2+5y^2+10z^2$,\\
$\left<1,3\right>$ in $\Qs{-1}$ &$\sim$& $w^2+x^2+3y^2+3z^2$,\\
$\left<1,3\right>$ in $\Qs{-2}$ &$\sim$& $w^2+3x^2+3y^2+6z^2$,\\
$\left<1,4\right>$ in $\Qs{-2}$ &$\sim$& $w^2+2x^2+4y^2+8z^2$,\\
$\left<1,5\right>$ in $\Qs{-2}$ &$\sim$& $w^2+2x^2+5y^2+10z^2$.
\end{tabular}\end{equation}
The universality of each of the forms on the right-hand side of (\ref{diags}) was shown by Ramanujan \cite{Ra}.  Thus, we have the universality of the Hermitian forms on the left-hand side of (\ref{diags}) immediately.

This leaves only eight other diagonal Hermitian lattices in (\ref{lat}),
\begin{equation}\label{non-diags}
\begin{tabular}{rcl}
$\left<1,1\right>$ in $\Qs{-3}$ &$\sim$& $w^2+wx+x^2+y^2+yz+z^2$,\\
$\left<1,1\right>$ in $\Qs{-7}$ &$\sim$& $w^2+wx+2x^2+y^2+yz+2z^2$,\\
$\left<1,1\right>$ in $\Qs{-11}$ &$\sim$& $w^2+wx+3x^2+y^2+yz+3z^2$,\\
$\left<1,2\right>$ in $\Qs{-3}$ &$\sim$& $w^2+wx+x^2+2y^2+2yz+2z^2$,\\
$\left<1,2\right>$ in $\Qs{-7}$ &$\sim$& $w^2+wx+2x^2+2y^2+2yz+4z^2$,\\
$\left<1,2\right>$ in $\Qs{-11}$ &$\sim$& $w^2+wx+3x^2+2y^2+2yz+6z^2$,\\
$\left<1,2\right>$ in $\Qs{-19}$ &$\sim$& $w^2+wx+5x^2+2y^2+2yz+10z^2$,\\
$\left<1,3\right>$ in $\Qs{-7}$ &$\sim$& $w^2+wx+2x^2+3y^2+3yz+6z^2$.\\
\end{tabular}\end{equation}
Of the quadratic forms on the right-hand side of \ref{non-diags}, only the last, $$w^2+wx+2x^2+3y^2+3yz+6z^2,$$ is an escalator lattice used in the proof of the 290-Theorem (see the table of quaternary escalators found in \cite{290}).  Thus, we may invoke the 290-Theorem to show the universality of the other seven quadratic forms in (\ref{non-diags}); the check that each of these forms represents all of $S_{290}$ is an easy computation.  It then follows directly that the first seven Hermitian forms in (\ref{non-diags}) are all universal.  

As the universality of $w^2+wx+2x^2+3y^2+3yz+6z^2$ is required in the proof of the 290-Theorem, we must still invoke an an \emph{ad hoc} argument to prove the universality of $\left<1,3\right>$ in $\Qs{-7}$.  One such argument is given in \cite{KP}.

Now, we turn to the non-diagonal Hermitian lattices in (\ref{lat}).  These are the remaining eight lattices,
\begin{equation}\label{non-non-diags}
\begin{tabular}{rcl}
$\left<1\right>\bot\left(\begin{array}{cc}2&-1+\omega_5\\-1+\bar{\omega}_5&3\end{array}\right)$ in $\Qs{-5}$ &$\sim$&$w^2+2x^2+2xy+3y^2+5z^2$,\\
$\left<1\right>\bot\left(\begin{array}{cc}2&\omega_6\\\bar{\omega}_6&3\end{array}\right)$ in $\Qs{-6}$ &$\sim$&$w^2+2x^2+3y^2+6z^2$,\\
$\left<1\right>\bot\left(\begin{array}{cc}2&\omega_{10}\\\bar{\omega}_{10}&5\end{array}\right)$ in $\Qs{-10}$ &$\sim$&$w^2+2x^2+3y^2+10z^2$,\\
$\left<1\right>\bot\left(\begin{array}{cc}2&\omega_{15}\\\bar{\omega}_{15}&2\end{array}\right)$ in $\Qs{-15}$ &$\sim$&$w^2+2x^2+xy+2y^2+wz+4z^2$,\\
$\left<1\right>\bot\left(\begin{array}{cc}2&\omega_{23}\\\bar{\omega}_{23}&3\end{array}\right)$ in $\Qs{-23}$ &$\sim$&$w^2+2x^2+xy+3y^2+wz+6z^2$,\\
$\left<1\right>\bot\left(\begin{array}{cc}2&-1+\omega_{23}\\-1+\bar{\omega}_{23}&3\end{array}\right)$ in $\Qs{-23}$ &$\sim$&$w^2+2x^2+xy+3y^2+wz+6z^2$,\\
$\left<1\right>\bot\left(\begin{array}{cc}2&\omega_{31}\\\bar{\omega}_{31}&4\end{array}\right)$ in $\Qs{-31}$ &$\sim$&$w^2+2x^2+xy+4y^2+wz+8z^2$,\\
$\left<1\right>\bot\left(\begin{array}{cc}2&-1+\omega_{31}\\-1+\bar{\omega}_{31}&4\end{array}\right)$ in $\Qs{-31}$ &$\sim$&$w^2+2x^2+xy+4y^2+wz+8z^2$.
\end{tabular}\end{equation}
Now, all of the diagonal quadratic forms in (\ref{non-non-diags}) are found in the list of universal forms obtained by Ramanujan \cite{Ra}.  Furthermore, none of the non-diagonal quadratic forms in (\ref{non-non-diags}) are escalator lattices, whence the universality of these forms follows directly from the 290-Theorem.  It then follows immediately that the Hermitian forms in (\ref{non-non-diags}) are all universal.
\end{proof}

\section{Remarks}
Kim, Kim, and Park \cite{KKP} have recently announced a criterion which completely characterizes the universality of Hermitian forms: \begin{15thmHerm}A positive-definite integral Hermitian form is universal if and only if it represents the ten integers $$\{1,2,3,5,6,7,10,13,14,15\}.$$
\end{15thmHerm}
Unfortunately, the proof of this result relies on our Main Theorem.  Consequently, Kim, Kim, and Park \cite{KKP}'s 15-Theorem for Hermitian Forms cannot be used to verify the universality of any of the forms in (\ref{lat}).

Kim, Kim, and Park \cite{KKP} do note that the 290-Theorem can be used to simplify some of the arguments in their proof of the 15-Theorem for Hermitian Forms.  Such simplifications would take the same form as those we have presented here to unify the classification of universal binary  Hermitian forms.

\end{document}